\documentclass[11pt]{article}

\usepackage[english]{babel} 
\usepackage[utf8]{inputenc} 
\usepackage[T1]{fontenc} 
\usepackage{fullpage}
\usepackage{amsfonts}
\usepackage{amsthm}
\usepackage{amsmath}
\usepackage{lmodern}
\usepackage{amssymb}
\usepackage{mathtools}
\usepackage{enumerate}
\usepackage{tikz-cd}
\usepackage{tikz}
\usepackage{graphicx}
\usepackage{stmaryrd}
\usepackage{longtable,graphics,multirow}
\usepackage{mathrsfs}
\usepackage{hyperref}
\usepackage{xcolor}

\hypersetup{
colorlinks=true,
linkcolor=[rgb]{0,0.4,0.6},%escala=255
filecolor=blue,
urlcolor=cyan,
citecolor=[rgb]{0,0.6,0},
anchorcolor=blue
}

\newtheorem{teor}{Theorem}
\numberwithin{teor}{section}

\newtheorem{lemma}[teor]{Lemma}

\newtheorem{prop}[teor]{Proposition}

\newtheorem{coro}[teor]{Corollary}

\newtheorem{remark}[teor]{Remark}

\theoremstyle{definition}\newtheorem{defi}[teor]{Definition}

\theoremstyle{definition}\newtheorem{example}{Example}
\numberwithin{example}{subsection}

\definecolor{customPink}{rgb}{1, 0.2, 0.4}%escala=255
\definecolor{customGreen}{rgb}{0, 0.6, 0}
\definecolor{customBlue}{rgb}{0, 0.4, 0.6}

\newcommand{\fim}{ $_\blacksquare$}
\def\P{\mathbb{P}}

\def\Z{\mathbb{Z}}

\def\Q{\mathbb{Q}}

\newcommand{\TypeA}{
\begin{tikzpicture}[thick,scale=1.4, every node/.style={scale=1}]

\draw[fill=white] (0,0) circle (1.5	pt);

\draw (0,0) node[above]{\small{$1$}};	

\draw (0.5,0) node{$\bullet$}; 

\draw (0.5,0) node[above]{\small{$1$}};	

\draw (1,0) node{$\bullet$};

\draw (1,0) node[above]{\small{$1$}};	

\draw (1.5,0) node{$\bullet$};

\draw (1.5,0) node[above]{\small{$1$}};

\draw (2,0) node{$\bullet$};

\draw (2,0) node[above]{\small{$1$}};

\draw[fill=white] (2.5,0) circle (1.5pt);

\draw (2.5,0) node[above]{\small{$1$}};

\draw (0.05,0) -- (0.5,0);

\draw (0.5,0) -- (1,0);

\draw[dashed] (1,0) -- (1.5,0);

\draw (1.5,0) -- (2,0);

\draw (2,0) -- (2.45,0);
\end{tikzpicture}
}

\newcommand{\TypeB}{
\begin{tikzpicture}[thick,scale=1.4, every node/.style={scale=1}]

\draw[fill=white] (0,0) circle (1.5	pt);

\draw (0,0) node[above]{\small{$2$}};	

\draw (0.5,0) node{$\bullet$}; 

\draw (0.5,0) node[above]{\small{$2$}};	

\draw (1,0) node{$\bullet$};

\draw (1,0) node[above]{\small{$2$}};	

\draw (1.5,0) node{$\bullet$};

\draw (1.5,0) node[above]{\small{$2$}};

\draw (2,0) node{$\bullet$};

\draw (2,0) node[above]{\small{$2$}};

\draw (2.5,0.3) node{$\bullet$};

\draw (2.5,0.3) node[above]{\small{$1$}};

\draw (2.5,-0.3) node{$\bullet$};

\draw (2.5,-0.3) node[below]{\small{$1$}};

\draw (0.05,0) -- (0.5,0);

\draw (0.5,0) -- (1,0);

\draw[dashed] (1,0) -- (1.5,0);

\draw (1.5,0) -- (2,0);

\draw (2,0) -- (2.5,0.3);

\draw (2,0) -- (2.5,-0.3);
\end{tikzpicture}
}

\newcommand{\TypeC}{
\begin{tikzpicture}[thick,scale=1.4, every node/.style={scale=1}]

\draw (0,0) node{$\bullet$};

\draw (0,0) node[above]{\small{$1$}};

\draw[fill=white] (0.5,0) circle (1.5pt);

\draw (0.5,0) node[above]{\small{$2$}};	

\draw (1,0) node{$\bullet$};

\draw (1,0) node[above]{\small{$1$}};	

\draw (0,0) -- (0.45,0);

\draw (0.55,0) -- (1,0);
\end{tikzpicture}
}

\newcommand{\TypeD}{
\begin{tikzpicture}[thick,scale=1.4, every node/.style={scale=1}]

\draw[fill=white] (0,0) circle (1.5pt);

\draw (0,0) node[above]{\small{$1$}};	

\draw[fill=white] (0.5,0) circle (1.5pt);

\draw (0.5,0) node[above]{\small{$1$}};	

\draw (0.05,0) -- (0.45,0);
\end{tikzpicture}
}

\newcommand{\TypeE}{
\begin{tikzpicture}[thick,scale=1.4, every node/.style={scale=1}]

\draw (0,0) node{$\star$};

\draw (0,0) node[above]{\small{$1$}};
\end{tikzpicture}
}

\newcommand{\DiagramX}{
\begin{tikzpicture}[thick,scale=1.4, every node/.style={scale=1}]

\draw[fill=white] (0,0) circle (1.5	pt);

\draw (0,0) node[above]{\small{$2$}};	

\draw (0.5,0) node{$\bullet$}; 

\draw (0.5,0) node[above]{\small{$n_1$}};	

\draw (1,0) node{$\bullet$};

\draw (1,0) node[above]{\small{$n_2$}};	

\draw (1.5,0) node{$\bullet$};

\draw (1.5,0) node[above]{\small{$n_m$}};	

\draw (0.05,0) -- (0.5,0);

\draw (0.5,0) -- (1,0);

\draw[dashed] (1,0) -- (1.5,0);
\end{tikzpicture}
}

\newcommand{\DiagramY}{
\begin{tikzpicture}[thick,scale=1.4, every node/.style={scale=1}]

\draw[fill=white] (0,0) circle (1.5	pt);

\draw (0,0) node[above]{\small{$2$}};	

\draw (0.5,0) node{$\bullet$}; 

\draw (0.5,0) node[above]{\small{$n_1$}};	

\draw (1,0) node{$\bullet$};

\draw (1,0) node[above]{\small{$n_2$}};	

\draw (1.5,0) node{$\bullet$};

\draw (1.5,0) node[above]{\small{$n_m$}};

\draw (2,0.3) node[right]{branch 1};

\draw (2,-0.3) node[right]{branch 2};

\draw (0.05,0) -- (0.5,0);

\draw (0.5,0) -- (1,0);

\draw[dashed] (1,0) -- (1.5,0);

\draw (1.5,0) -- (2.06,0.3);

\draw (1.5,0) -- (2.06,-0.3);
\end{tikzpicture}
}

\newcommand{\DiagramZ}{
\begin{tikzpicture}[thick,scale=1.4, every node/.style={scale=1}]

\draw[fill=white] (0,0) circle (1.5	pt);

\draw (0,0) node[above]{\small{$2$}};	

\draw (0.5,0) node{$\bullet$}; 

\draw (0.5,0) node[above]{\small{$n_1$}};	

\draw (1,0.4) node{$\bullet$};	

\draw (1,0.4) node[above]{\small{$n_2$}};	

\draw (1.5,0.2) node{$\bullet$};

\draw (1.5,0.2) node[above]{\small{$n_3$}};	

\draw (1.5,-0.2) node{$\bullet$};

\draw (1.5,-0.2) node[below]{\small{$n_4$}};	

\draw (1,-0.4) node{$\bullet$};

\draw (1,-0.4) node[below]{\small{$n_m$}};	

\draw (0.05,0)--(0.5,0);

\draw (0.5,0)--(1,0.4);

\draw (1,0.4)--(1.5,0.2);

\draw (1.5,0.2)--(1.5,-0.2);

\draw[dashed] (1.5,-0.2)--(1,-0.4);

\draw (1,-0.4)--(0.5,0);
\end{tikzpicture}
}

\begin{document}

\author{Renato Dias Costa}
\date{}
\title{Classification of Conic Bundles on a Rational\\Elliptic Surface in any Characteristic}
\maketitle

\begin{abstract}
Let $X$ be a rational elliptic surface with elliptic fibration $\pi:X\to\P^1$ over an algebraically closed field $k$ of any characteristic. Given a conic bundle $\varphi:X\to\P^1$ we use numerical arguments to classify all possible fibers of $\varphi$ and 	study the interplay between singular fibers of $\pi$ and $\varphi$.
\end{abstract}

\tableofcontents
\newpage
\section{Introduction}\
\indent Let $S$ be a smooth, projective surface over a field $k$. We define a conic bundle on $S$ in the most natural way, namely as a morphism onto a curve $S\to C$ whose general fiber is a smooth, rational curve. A prominent case where conic bundles arise is in the classification of minimal models of $k$-rational surfaces. Iskovskikh showed \cite[Theorem 1]{Isko} that if $S$ is a $k$-minimal rational surface, then $S$ is either (a) $\P^2_k$, (b) a quadric on $\P^3_k$, (c) a del Pezzo surface, or (d) $S$ admits a conic bundle such that its singular fibers are isomorphic to a pair of lines meeting at a point. The latter was named a \textit{standard conic bundle} by Manin and Tsfasman \cite[Subsection 2.2]{ManTsfa}. We remark that the notion of standard conic bundle is often extended to higher dimension and plays a role in the classification of threefolds over $\Bbb{C}$ in the Minimal Model Program (see, for example, \cite{Sarkisov}, \cite{IskoRationalityProblem} and the survey \cite{Prokho}).

We are concerned more specifically with conic bundles on a rational surface $X$ over an algebraically closed field $k$ of any characteristic with an elliptic fibration $\pi:X\to\P^1$, i.e. a genus $1$ fibration with a section. This is motivated by results and techniques from \cite{Salgado}, \cite{LoughSalgado}, \cite{AliceSalgado}, \cite{AliceSalgadoII} and \cite{Alice}, which we now briefly explain.

In \cite{Salgado}, $k$ is a number field and the goal is to study sets of fibers $\pi^{-1}(t)$ with $t\in\P^1(k)$ such that the Mordell-Weil rank $r_t$ of the fiber is greater than the generic rank $r$ of the fibration. The presence of a bisection, i.e. a rational curve $C\subset X$ such that $C\cdot(-K_X)=2$ induces a family of bisections, which turns out to be a conic bundle on $X$. The author explores the existence of two such conic bundles and proves that the sets $\{r_t\geq r+1\}$ and $\{r_t\geq r+2\}$ are infinite under certain hypothesis. This strategy was later refined in \cite{LoughSalgado} --- using, for example, one conic bundle instead of two --- in order to study the structure of these sets in view of the Hilbert property \cite[Chapter 3]{Serre}. Further developments are found in \cite{DiasCostaSalgado}, where special conic bundles are considered in order to study the set $\{r_t\geq r+3\}$, and also in \cite{HindrySalgado}, where many ideas are extended from elliptic fibrations to families of Abelian varieties over a rational curve. 

When $k$ is algebraically closed with characteristic zero, the existence of conic bundles on $X$ is also used in \cite{AliceSalgado}, \cite{AliceSalgadoII} in order to classify elliptic fibrations on K3 surfaces which are quadratic covers of $X$. More precisely, given a degree two morphism $f:\P^1\to\P^1$ ramified away from nonreduced fibers of $\pi$, the induced K3 surface is $X':=X\times_f\P^1$. The base change also gives rise to an elliptic fibration $\pi':X'\to\P^1$ and a degree two map $f':X'\to X$. By composition with $f'$, every conic bundle on $X$ induces a genus $1$ fibration on $X'$, oftentimes having a section, in which case we get elliptic fibrations distinct from $\pi$.

In a different context, conic bundles also appear in \cite{Alice}. Here $k=\Bbb{C}$ and the goal is to find generators for the Cox ring $\mathcal{R}(X):=\bigoplus_{[D]}H^0(X,\mathcal{O}_X(D))$, where $[D]$ runs through $\text{Pic}(X)$. Given a rational elliptic fibration $\pi:X\to\P^1$, the ring $\mathcal{R}(X)$ is finitely generated if and only if $X$ is a Mori Dream Space \cite[Proposition 2.9]{HuKeel}, which in turn is equivalent to $\pi$ having generic rank zero \cite[Corollary 5.4]{ArtebaniLaface}. Assuming this is the case, the authors show that in many configurations of $\pi$ each minimal set of generators of $\mathcal{R}(X)$ must contain an element $g\in H^0(X,\mathcal{O}_X(D))$, where $D$ is a fiber of a conic bundle on $X$, whose possibilities are explicitly described.

This argues for a detailed analysis of conic bundles on rational elliptic surfaces, which is what we aim for. We explore ideas that are either implicit or mentioned in passing throughout \cite{Salgado}, \cite{LoughSalgado}, \cite{AliceSalgado}, \cite{AliceSalgadoII}, \cite{Alice} and make generalizations whenever they apply. 

The paper is organized as follows. We set the preliminaries in Section~\ref{setup}, where we define an elliptic surface, state some of its properties, define conic bundles and list some general facts about divisors, fibers and linear systems. Section~\ref{NumChar} is dedicated to showing that conic bundles can be characterized numerically through a bijective correspondence with certain classes of the Néron-Severi group, in a sense made precise in Theorem~\ref{correspondence}. In Section~\ref{Classification} we prove our first result, which is the complete classification of conic bundle fibers in arbitrary characteristic.
\\ \\
\noindent\textbf{Theorem \color{black}\ref{teor_classif}.} \textit{Let $X$ be a rational elliptic surface with elliptic fibration $\pi:X\to\P^1$ and let $\varphi:X\to\P^1$ be a conic bundle. If $D$ is a fiber of $\varphi$, then the intersection graph of $D$ fits one of the following types. Conversely, if $D$ fits any of these types, then $|D|$ induces a conic bundle $X\to\P^1$.}
\begin{table}[h]
\begin{center}
\centering
\begin{tabular}{|c|c|} 
\hline
Type & Intersection Graph\\ 
\hline
\multirow{3}{*}{\hfil 0} & \multirow{3}{*}{\hfil \TypeE}\\ 
&\\
&\\
\hline
\multirow{3}{*}{\hfil $A_2$} & \multirow{3}{*}{\hfil \TypeD}\\
& \\
& \\
\hline
\multirow{3}{*}{\hfil $A_n\,(n\geq 3)$} & \multirow{3}{*}{\hfil \TypeA}\\
& \\ 
& \\
\hline
\multirow{3}{*}{\hfil $D_3$} & \multirow{3}{*}{\hfil \TypeC}\\
& \\
& \\
\hline
\multirow{4}{*}{\hfil $D_m\,(m\geq 4)$} & \multirow{4}{*}{\hfil \TypeB}\\ %& \multirow{4}{*}{\hfil $2P+\sum_{i=0}^{n-2}2\Theta_i+\Theta_{n-1}+\Theta_n$}\\
& \\ 
& \\
& \\
\hline
\end{tabular}
\end{center}
\begin{align*}
\star&\,\,\text{smooth, irreducible curve of genus zero}\\
\circ&\,(-1)\text{-curve (section of }\pi)\\
\bullet&\,(-2)\text{-curve (component of a reducible fiber of }\pi)
\end{align*}
\end{table}

In Section~\ref{conic_vs_elliptic} we investigate how the fibers of $\pi$ affect the possibilities for conic bundles on $X$ and obtain our second result.
\\ \\
\noindent\textbf{Theorem \ref{fiber_table}.}\textit{Let $X$ be a rational elliptic surface with elliptic fibration $\pi:X\to\P^1$. Then the following statements hold:}
\begin{enumerate}[a)]
\item $X$ admits a conic bundle with an $A_2$ fiber $\Rightarrow $ $\pi$ is not extremal (i.e. has positive generic rank).
\item $X$ admits a conic bundle with an $A_n$ $(n\geq 3)$ fiber $\Leftrightarrow$ $\pi$ has a reducible fiber other than {\normalfont $\text{II}^*$}.
\item $X$ admits a conic bundle with a $D_3$ fiber $\Leftrightarrow$ $\pi$ has at least two reducible fibers.	
\item $X$ admits a conic bundle with a $D_m$ $(m\geq 4)$ fiber $\Leftrightarrow$ $\pi$ has a nonreduced fiber or a fiber {\normalfont $\text{I}_n$ $(n\geq 4)$}.
\end{enumerate}

At last, in Section~\ref{construction}, we describe a method for constructing conic bundles from pencils of plane curves with genus $0$, which we apply in Section~\ref{examples} to produce examples.
\\ \\
\noindent\textbf{Acknowledgements.} The author was supported by FAPERJ (grant E-26/201.182/2020). He thanks Cecília Salgado for suggesting the problem and for her guidance along the writing process; Alice Garbagnati, for the generosity of reading this work carefully and proposing follow-up ideas.

\section{Set up}\label{setup}\	
\indent In what follows, all surfaces are projective and smooth over an algebraically closed field $k$ of any characteristic, unless otherwise specified. We always use $X$ to denote an \hyperref[ellipticsurfaces]{elliptic surface}. In Subsection~\ref{EllipticSurfaces} we define elliptic surfaces and refer to classical results. In Subsection~\ref{RES} we cover some well-known facts about rational elliptic surfaces and present two additional lemmas. Conic bundles are defined in Subsection~\ref{ConicBundles}.

\subsection{Elliptic surfaces}\label{EllipticSurfaces}
\begin{defi}\label{ellipticsurfaces}
A surface $X$ is called an \textit{elliptic surface} if there is smooth projective curve $C$ and a surjective morphism $\pi:X\to C$, called an \textit{elliptic fibration}, such that
\begin{enumerate}[i)]
\item The fiber $\pi^{-1}(t)$ is a smooth genus $1$ curve for all but finitely many $t\in C$.
\item (existence of a section) There is a morphism $\sigma:C\to X$ such that $\pi\circ\sigma=\text{id}_C$, called a \textit{section}.
\item (relative minimality)  No fiber of $\pi$ contains an exceptional curve in its support (i.e., a smooth rational curve with self-intersection $-1$).\label{rel_min}
\end{enumerate}
\end{defi}

\begin{remark}\label{relative_minimality}
\normalfont Condition \ref{rel_min}) can be understood as an extra hypothesis for $\pi$. For our purposes it is a natural one, since it assures that the fibers are as in Kodaira's classification (Theorem \ref{Kodaira_fibers}) and that some convenient properties hold for rational elliptic surfaces (see Theorem~\ref{RESpencil} and Theorem~\ref{anticanonical}).
\end{remark}
\begin{remark}\label{Sections}
\normalfont Given $\sigma:C\to X$ as in Definition~\ref{ellipticsurfaces}, the curve $P:=\sigma(C)$ on $X$ is also called a section. The image $\sigma(C)$ is isomorphic to $C$ and meets a general fiber of $\pi$ at one point. Conversely, for every smooth curve $P\subset X$ meeting a general fiber of $\pi$ at one point, there is a morphism $\sigma:C\to X$ such that $\pi\circ \sigma=\text{id}_C$ and $\sigma(C)=P$ \cite[Subsection 3.4]{Schuett-Shioda}.
\end{remark}
\newpage
We refer to some classical results on elliptic surfaces.

\begin{teor}{\normalfont \cite[Section 6]{Tate}}\label{Kodaira_fibers}
Let $\pi:X\to C$ be an elliptic fibration. If $F$ a singular fiber of $\pi$, then all possibilities for $F$ are listed below. The symbols $\widetilde{A}_n, \widetilde{D}_n, \widetilde{E}_n$ indicate the type of extended Dynking diagram formed by the intersection graph of $F$.
\begin{center}
\includegraphics[scale=0.65]{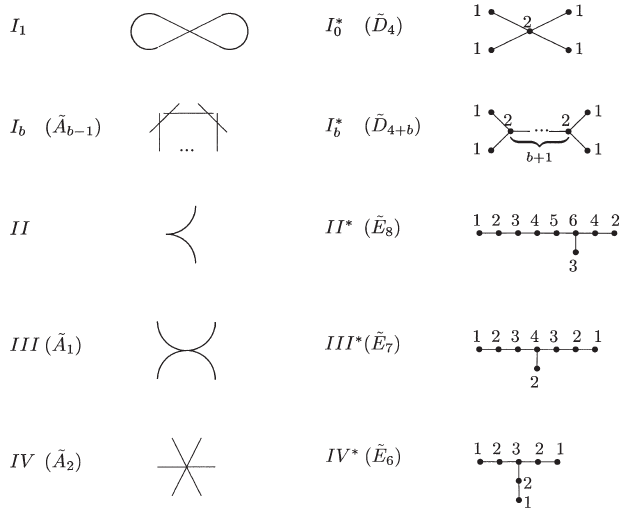}
\end{center}
Moreover, when $F$ is irreducible (namely of types $\text{I}_1$ and $\text{II}$), it is a singular, integral curve of arithmetic genus $1$. When $F$ is reducible, all members in its support are smooth, rational curves with self-intersection $-2$.
\end{teor}

\begin{remark}\label{remarkchar}
\normalfont Tate's algotithm for identifying singular fibers involves local Weierstrass forms and is valid over perfect fields (e.g. algebraically closed). This local analysis can be pathological in characteristics $2$, $3$ but even then all possible singular fibers are still covered by the list in Theorem~\ref{Kodaira_fibers} \cite[Section 6]{Tate}. We note that for non-perfect fields, new fiber types may occur (see \cite{Szydlo}). In this paper we do not deal with local Weierstrass forms but only with the numerical behavior of fibers as divisors, so we are free to use Theorem~\ref{Kodaira_fibers} over algebraically closed fields of any characteristic.
\end{remark}

\begin{teor}{\normalfont \cite[Thm. 3.1]{Shioda}}\label{numalgequiv}
On an elliptic surface $X$, algebraic and numerical equivalences coincide, i.e. {\normalfont $D_1,D_2\in\text{Div}(X)$} are equivalent in {\normalfont $\text{NS}(X)$} $\Leftrightarrow D_1\cdot D=D_2\cdot D$ for every {\normalfont $D\in\text{Div}(X)$}.
\end{teor}

\newpage

\subsection{Rational elliptic surfaces}\label{RES}
\begin{defi}
We say that $X$ is a \textit{rational elliptic surface} if there is an elliptic fibration $\pi:X\to C$ and $X$ is a rational surface.
\end{defi}

We exhibit a standard method for constructing rational elliptic surfaces. Since our base field $k$ is algebraically closed, every elliptic surface over $k$ can be obtained using this method (Theorem~\ref{RESpencil}), which is not the case over arbitrary fields. 

The construction works as follows. Let $F,G$ be cubics on $\Bbb{P}_k^2$, at least one of them smooth. The intersection $F\cap G$ has nine points counted with multiplicity and the pencil of cubics $\mathcal{P}:=\{sF+tG=0\mid (s:t)\in\P^1\}$ has $F\cap G$ as base locus. Let $\phi:\P^2\dasharrow\P^1$ be the rational map associated to $\mathcal{P}$. The blowup $p:X\to \Bbb{P}^2$ at the base locus resolves the indeterminacies of $\phi$ and yields a surface with an elliptic fibration $\pi:X\to\Bbb{P}^1$.

\begin{displaymath}
\begin{tikzcd}
X\arrow{r}{p}\arrow[bend right, swap]{rr}{\pi} & \P^2\arrow[dashed]{r}{\phi} &\P^1
\end{tikzcd}
\end{displaymath}

By construction, $X$ is rational. Conversely, every rational elliptic surface over $k$ (which is algebraically closed) can be obtained by this method.

\begin{teor}{\normalfont \cite[Theorem 5.6.1]{Dolga}}\label{RESpencil}
Every rational elliptic surface over an algebraically closed field is isomorphic to the blowup of $\P^2$ at the base locus of a pencil of cubics.
\end{teor}

We mention some other distinguished features of rational elliptic surfaces.

\begin{teor}{\normalfont \cite[Section 8.2]{Schuett-Shioda}}\label{anticanonical}
Let $X$ be a rational elliptic surface and $\pi:X\to\P^1$ its elliptic fibration. Then
\begin{enumerate}[i)]
\item $\chi(X,\mathcal{O}_X)=1$, where $\chi(X,\mathcal{O}_X)=h^0(X,\mathcal{O}_X)-h^1(X,\mathcal{O}_X)+h^2(X,\mathcal{O}_X)$.
\item $-K_X$ is linearly equivalent to any fiber of $\pi$. In particular, $-K_X$ is nef.
\item Every section of $\pi$ is an exceptional curve (smooth, rational curve with self-intersection $-1$).
\end{enumerate}
\end{teor}

We include two elementary results. Lemma \ref{detect_fibers} provides a simple test to detect fibers of the elliptic fibration and Lemma \ref{neg_curves} describes the negative curves on a rational elliptic surface.

\begin{lemma}\label{detect_fibers}
Let $\pi:X\to\P^1$ be a rational elliptic fibration and $E$ an integral curve in $X$. If $E\cdot K_X=0$, then $E$ is a component of a fiber of $\pi$. If moreover $E^2=0$, then $E$ is a fiber.
\end{lemma}
\noindent\textit{Proof.} If $P\in E$, then the fiber $F:=\pi^{-1}(\pi(P))$ intersects $E$ at $P$, i.e. $E\cap F\neq\emptyset$. On the other hand, $-K_X$ is linearly equivalent to $F$ by Theorem \ref{anticanonical}, so $E\cdot F=-E\cdot K_X=0$. Hence $E$ must be a component of $F$. Assuming moreover that $E^2=0$, we prove that $E=F$. In case $F$ is smooth, this is clear. So we assume $F$ is singular and analyze its Kodaira type according to Theorem~\ref{Kodaira_fibers}. Notice that $F$ is not reducible, otherwise $E$ would be a $(-2)$-curve, which contradicts $E^2=0$. Hence $F$ is either of type $\text{I}_1$ or $\text{II}$. In both cases $F$ is an integral curve, therefore $E=F$.\fim
\newpage
\begin{lemma}\label{neg_curves}
Let $\pi:X\to\P^1$ be a rational elliptic fibration. Every negative curve on $X$ is either a $(-1)$-curve (section of $\pi$) or a $(-2)$-curve (component from a reducible fiber of $\pi$). 
\end{lemma}
\noindent\textit{Proof.} Let $E$ be any integral curve in $X$ with $E^2<0$. By Theorem \ref{anticanonical}, $-K_X$ is nef and linearly equivalent to any fiber of $\pi$. So $E\cdot(-K_X)\geq 0$ and by adjunction $2p_a(E)-2=E^2+E\cdot K_X<0$, which only happens if $p_a(E)=0$. Consequently $E^2=-1$ or $-2$. In case $E^2=-2$ we have $E\cdot K_X=0$, so $E$ is fiber a component by Lemma \ref{detect_fibers} and this fiber is reducible by Theorem \ref{Kodaira_fibers}. If $E^2=-1$, by adjunction $E\cdot(-K_X)=1$. But $-K_X$ is lineary equivalent to any fiber, so $E$ meets a general fiber at one point, therefore $E$ is a section by Remark \ref{Sections}.\fim

\subsection{Conic bundles}\label{ConicBundles}\
\indent The following definition of conic bundle is essentially identical to the ones in \cite[Definition~2.5]{LoughSalgado}, \cite[Definition 3.1]{AliceSalgado}. The difference here is that $k$ is algebraically closed with arbitrary characteristic.

\begin{defi}\label{def_conicbundles}
A \textit{conic bundle} on a surface $S$ is a surjective morphism onto a smooth curve $\varphi:S\to C$ whose general fiber is a smooth, irreducible curve of genus $0$.
\end{defi}

\begin{remark}
\normalfont When $S$ is rational, $C$ is isomorphic to $\P^1$ by Lüroth's theorem. 
\end{remark}

\begin{remark}
\normalfont Definition \ref{def_conicbundles} is more general than the one of standard conic bundle in \cite[Subsection 2.2]{ManTsfa}, in which every singular fiber is isomorphic to a pair of lines meeting at a point. Indeed we show in Theorem~\ref{teor_classif} that such a singular fiber is one among four possible types using our definition.
\end{remark}

\begin{remark}
\normalfont Definition \ref{def_conicbundles} is geometric in nature. In Section \ref{NumChar} we show that a conic bundle may be equivalently defined numerically by a certain Néron-Severi class $[D]\in \text{NS}(X)$.
\end{remark}
\subsection{General facts about divisors, fibers and linear systems}\
\indent In this subsection we let $S$ be a smooth projective surface (not necessarily elliptic) over an algebraically closed field. We list some elementary results involving divisors, fibers of morphisms $S\to\P^1$ and linear systems, none of which depends on the theory of elliptic surfaces.

\begin{lemma}\label{critical_nef}
Let $D$ be an effective divisor such that $D\cdot C=0$ for all {\normalfont $C\in\text{Supp }D$}. Then $D$ is nef and $D^2=0$.
\end{lemma}
\noindent\textit{Proof.} To see that $D$ is nef, just notice that $D\cdot C'\geq 0$ for every $C'\notin\text{Supp }D$. To prove that $D^2=0$ let $D=\sum_in_iC_i$, where each $C_i$ is in $\text{Supp }D$. Then $D^2=D\cdot(\sum_in_iC_i)=\sum_in_iD\cdot C_i=0$.\fim
\\ \\
\indent The next lemma explores some properties of fibers of a surjective morphism $f:S\to\P^1$.

\begin{lemma}\label{lemma_fibers}
Let $F$ be a fiber of a surjective morphism $f:S\to\P^1$. Then the following hold:
\begin{enumerate}[a)]
\item $F\cdot C=0$ for all {\normalfont $C\in\text{Supp }F$}.
\item If $F$ is connected and $E$ is a divisor such that {\normalfont $\text{Supp }E\subset\text{Supp }F$}, then $E^2\leq 0$ and $E^2=0$ if and only if $E=rF$ for some $r\in\Q$.
\item If $F_1,...,F_n$ are the connected components of $F$, then each $F_i$ is nef with $F_i^2=0$ and $F_i\cdot K_S\in 2\Z$.
\end{enumerate}
\end{lemma}

\noindent\textit{Proof.} a) Taking an arbitrary $C\in\text{Supp }F$ and another fiber $F'\neq F$, we have $F\cdot C=F'\cdot C=0$.
\\ \\
\indent b) This is Zariski's lemma \cite[Ch. 6, Lemma 6]{Peters}.
\\ \\
\indent c) Fix $i$. To prove that $F_i$ is nef and $F_i^2=0$ we show that $F_i\cdot C_i=0$ for any $C_i\in\text{Supp }F_i$ then apply Lemma~\ref{critical_nef}. Indeed, if $j\neq i$, the components $F_i,F_j$ are disjoint, so $F_j\cdot C_i=0$. Since $F$ is a fiber and $C_i\in\text{Supp }F_i\subset\text{Supp }F$, then $F\cdot C_i=0$ by~a). Hence $F_i\cdot C_i=(F_1+...+F_n)\cdot C_i=F\cdot C_i=0$, as desired. For the last part, by Riemann-Roch $F_i\cdot K_S=2(\chi(S,\mathcal{O}_S)-\chi(S,\mathcal{O}_S(F_i)))\in 2\Z$.\fim
\\ \\
\indent This last lemma is a property of linear systems with no fixed components.
\begin{lemma}\label{fixed_component}
Let $E,E'$ be effective divisors such that $E'\leq E$ and that the linear systems $|E|,|E'|$ have the same dimension. If $|E|$ has no fixed components, then $E'=E$.
\end{lemma}
\noindent\textit{Proof.} The fact that $E'\leq E$ implies that $H^0(S,\mathcal{O}_S(E'))$ is a subspace of $H^0(S,\mathcal{O}_S(E))$. By hypothesis these spaces	 have the same dimension, so $H^0(S,\mathcal{O}_S(E'))=H^0(S,\mathcal{O}_S(E))$. Hence
\begin{align*}
|E|&=\{E+\text{div}(f)\mid f\in H^0(S,\mathcal{O}_S(E))\}\\
&=\{E'+\text{div}(f)+(E-E')\mid f\in H^0(S,\mathcal{O}_S(E'))\}\\
&=|E'|+(E-E').
\end{align*}

Assuming $|E|$ has no fixed components, we must have $E-E'=0$.\fim

\section{Numerical characterization of conic bundles}\label{NumChar}\
\indent Let $\pi:X\to\P^1$ be a rational elliptic fibration. We give a characterization of conic bundles on $X$ which shows the numerical nature of conic bundles on rational elliptic surfaces. The motivation for this comes from the following. Let $\varphi:X\to\P^1$ be a conic bundle and $C$ a general fiber of $\varphi$, which is a smooth, irreducible curve of genus zero. Clearly $C$ is a nef divisor with $C^2=0$ and by adjunction $C\cdot(-K_X)=2$. These three numerical properties are enough to prove that $|C|$ is a base point free pencil and consequently the induced morphism $\varphi_{|C|}:X\to\P^1$ is precisely $\varphi$.

Conversely, let $D$ be a nef divisor with $D^2=0$ and $D\cdot(-K_X)=2$. Since numerical and algebraic equivalence coincide by Theorem \ref{numalgequiv}, it makes sense to consider the class $[D]\in\text{NS}(X)$. The natural question is whether $[D]$ induces a conic bundle on $X$. The answer is yes, moreover there is a natural correspondence between such classes and conic bundles (Theorem~\ref{correspondence}), which is the central result of this section.

In order to prove this correspondence we need a numerical analysis of a given class $[D]\in\text{NS}(X)$ so that we can deduce geometric properties of the induced morphism $\varphi_{|D|}:X\to\P^1$, such as connectivity of fibers (Proposition \ref{connectedness}) and composition of their support (Proposition \ref{components_conicbundleclass}). These properties are also crucial to the classification of fibers in Section~\ref{Classification}.

\begin{defi}
A class $[D]\in\text{NS}(X)$ is called a \textit{conic class} when
\begin{enumerate}[i)]
\item $D$ is nef.
\item $D^2=0$.
\item $D\cdot (-K_X)=2$.
\end{enumerate}
\end{defi}

\begin{lemma}\label{Riemann-Roch}
Let {\normalfont $[D]\in\text{NS}(X)$} be a conic class. Then $|D|$ is a base point free pencil and therefore induces a surjective morphism $\varphi_{|D|}:X\to\P^1$.
\end{lemma}
\noindent\textit{Proof.} By \cite[Theorem III.1(a)]{Har}, $|D|$ is base point free and $h^1(X,D)=0$. We have $\chi(X,\mathcal{O}_X)=1$ by Theorem \ref{anticanonical}, and Riemann-Roch gives $h^0(X,D)+h^2(X,D)=2$, so we only need to prove $h^2(X,D)=0$. Assume by contradiction that $h^2(X,D)\geq 1$. By Serre duality $h^0(K_X-D)\geq 1$, so $K_X-D$ is linearly equivalent to an effective divisor. Since $D$ is nef, $(K_X-D)\cdot D\geq 0$, which contradicts $D^2=0$ and $D\cdot(-K_X)=2$.\fim

\begin{remark}
\normalfont It follows from Lemma~\ref{Riemann-Roch} that a conic class has an effective representative.
\end{remark}

Notice that we do not know a priori that the morphism $\varphi_{|D|}:X\to\P^1$ in Lemma~\ref{Riemann-Roch} is a conic bundle. At this point we can only say that if $C$ is a smooth, irreducible fiber of $\varphi_{|D|}$, then $g(C)=0$ by adjunction. However it is still not clear whether a general fiber of $\varphi_{|D|}$ is irreducible and smooth. We prove that this is the case in Proposition~\ref{finitely_many_sinuglar_fibers}. In order to do that we need information about the components of $D$ from Proposition~\ref{components_conicbundleclass} and the fact that $D$ is connected from Proposition~ \ref{connectedness}.

\begin{prop}\label{components_conicbundleclass}
Let {\normalfont $[D]\in\text{NS}(X)$} be a conic class. If $D$ is an effective representative, then every curve {\normalfont $E\in\text{Supp }D$} is a smooth rational curve with $E^2\leq 0$.
\end{prop}

\noindent\textit{Proof.} Take an arbitrary $E\in\text{Supp }D$. By Lemma~\ref{Riemann-Roch}, $D$ is a fiber of the morphism $\varphi_{|D|}:X\to\P^1$ induced by $|D|$, hence $D\cdot E=0$ by Lemma~\ref{critical_nef}. Assuming $E^2>0$ by contradiction, the fact that $D^2=0$ implies that $D$ is numerically equivalent to zero by the Hodge index theorem \cite[Thm. V.1.9, Rmk. 1.9.1]{Hartshorne}. This is absurd because $D\cdot(-K_X)=2\neq 0$, so indeed $E^2\leq 0$.

To show that $E$ is a smooth rational curve, it suffices to prove that $p_a(E)=0$. By Theorem \ref{anticanonical}, $-K_X$ is linearly equivalent to any fiber of $\pi$, in particular $-K_X$ is nef and $E\cdot K_X\leq 0$. By adjunction $2p_a(E)-2=E^2+E\cdot K_X\leq 0$, so $p_a(E)\leq 1$. Assume by contradition that $p_a(E)=1$. This can only happen if $E^2=E\cdot K_X=0$, so $E$ is a fiber of $\pi$ by Lemma \ref{detect_fibers}. In this case $E$ is linearly equivalent to $-K_X$, so $D\cdot E=D\cdot(-K_X)=2$, which contradicts $D\cdot E=0$.\fim
\\ \\
\indent While Proposition~\ref{components_conicbundleclass} provides information about the support of $D$, the next proposition states that $D$ is connected, which is an important fact about the composition of $D$ as a whole. 
\begin{prop}\label{connectedness}
Let {\normalfont $[D]\in\text{NS}(X)$} be a conic class. If $D$ is an effective representative, then $D$ is connected.
\end{prop}
\noindent\textit{Proof.} Let $D=D_1+...+D_n$, where $D_1,...,D_n$ are connected components. By Lemma~\ref{lemma_fibers} c) each $D_i$ is nef with $D_i^2=0$ and $D_i\cdot K_X\in 2\Z$. Since $-K_X$ is nef by Theorem~\ref{anticanonical} and $D\cdot(-K_X)=2$, then $D_{i_0}\cdot(-K_X)=2$ for some $i_0$ and $D_i\cdot (-K_X)=0$ for $i\neq i_0$. In particular $[D_{i_0}]\in\text{NS}(X)$ is a conic class. By Lemma~\ref{fixed_component}, both $|D|,|D_{i_0}|$ are pencils, so $D=D_{i_0}$ by Lemma~\ref{fixed_component}.\fim
\\ \\
\indent We use Propositions \ref{components_conicbundleclass} and \ref{connectedness} to conclude that $\varphi_{|D|}:X\to\P^1$ is indeed a conic bundle.
\newpage
\begin{prop}\label{finitely_many_sinuglar_fibers}
Let {\normalfont $[D]\in\text{NS}(X)$} be a conic class. Then all fibers of $\varphi_{|D|}:X\to\P^1$ are smooth, irreducible curves of genus $0$ except for finitely many which are reducible and supported on negative curves. In particular, $\varphi_{|D|}$ is a conic bundle.
\end{prop}
\noindent\textit{Proof.} Let $F$ a fiber of $\varphi_{|D|}$. Since $F$ is linearly equivalent to $D$, then $[F]=[D]\in\text{NS}(X)$, so $F$ is connected by Proposition \ref{connectedness}. By Proposition \ref{components_conicbundleclass} every $E\in\text{Supp }F$ has $g(E)=0$ and $E^2\leq 0$.

First assume $E^2=0$ for some $E\in\text{Supp }F$. Since $F$ is a connected fiber of $\varphi_{|D|}$, then $E=rF$ for some $r\in\Q$ by Lemma~\ref{lemma_fibers}. Because $F\cdot(-K_X)=2$, we have $E\cdot K_X=-2r$. By adjunction, $r=1$, so $F=E$ is a smooth, irreducible curve of genus $0$. 

Now assume $E^2<0$ for every $E\in\text{Supp }F$. Then $F$ must be reducible, otherwise $F^2=E^2<0$, which is absurd since $F$ is a fiber of $\varphi_{|D|}$. Conversely, if $F$ is reducible, then $E^2<0$ for all $E\in\text{Supp }F$, otherwise $E^2=0$ for some $E$ and by the last paragraph $F$ is irreducible, which is a contradiction. 

This shows that either $F$ is smooth, irreducible of genus $0$ or $F$ is reducible, in which case $F$ is supported on negative curves. We are left to show that $\varphi_{|D|}$ has finitely many reducible fibers.

Assume by contradiction that there is a sequence $\{F_n\}_{n\in\Bbb{N}}$ of distinct reducible fibers of $\varphi_{|D|}$. In particular each $F_n$ is supported on negative curves, which are either $(-1)$-curves (sections of $\pi$) or $(-2)$-curves (components of reducible fibers) by Theorem \ref{neg_curves}.

Since $\pi$ has finitely many singular fibers, the number of $(-2)$-curves in $X$ is finite, so there are finitely many $F_n$ with $(-2)$-curves in its support. Excluding such $F_n$, we may assume all members in $\{F_n\}_{n\in\Bbb{N}}$ are supported on $(-1)$-curves. For each $n$, take $P_n\in\text{Supp }F_n$. The fibers $F_n,F_m$ are disjoint, so $P_n,P_m$ are disjoint when $n\neq m$. If we successively contract the exceptional curves $P_1,...,P_\ell$, we are still left with an infinite set $\{P_n\}_{n>\ell}$ of exceptional curves, so we cannot reach a minimal model for $X$, which is absurd.\fim

\begin{remark}
\normalfont In characteristic zero the proof of Proposition \ref{finitely_many_sinuglar_fibers} can be made simpler by applying Bertini's theorem, which guarantees that the general fiber is smooth from the fact that $X$ is smooth.
\end{remark}

We now prove the numerical characterization of conic bundles.

\begin{teor}\label{correspondence}
Let $\pi:X\to\P^1$ be a rational elliptic fibration. If {\normalfont $[D]\in\text{NS}(X)$} is a conic class, then $|D|$ is a base point free pencil whose induced morphism $\varphi_{|D|}:X\to\P^1$ is a conic bundle. Moreover, the map $[D]\mapsto \varphi_{|D|}$ has an inverse $\varphi\mapsto [F]$, where $F$ is any fiber of $\varphi$. This gives a natural correspondence between conic classes and conic bundles.
\end{teor}

\textit{Proof.} Given a conic class $[D]\in\text{NS}(X)$, by Proposition \ref{finitely_many_sinuglar_fibers} the general fiber of $\varphi_{|D|}:X\to\P^1$ is a smooth, irreducible curve of genus $0$, so $\varphi_{|D|}$ is a conic bundle.

Conversely, if $\varphi:X\to\P^1$ is a conic bundle and $C$ is a smooth, irreducible fiber of $\varphi$, in particular $C^2=0$, $g(C)=0$ and by adjunction $C\cdot(-K_X)=2$. Clearly $C$ is nef, so $[C]\in\text{NS}(X)$ is a conic class. Moreover, any other fiber $F$ of $\varphi$ is linearly equivalent to $C$, therefore $[F]=[C]\in\text{NS}(X)$ and the map $\varphi\mapsto [F]$ is well defined.

We verify that the maps are mutually inverse. Given a class $[D]$ we may assume $D$ is effective since $|D|$ is a pencil, so $D$ is a fiber of $\varphi_{|D|}$ and $\varphi_{|D|}$ is sent back to $[D]$. Conversely, given a conic bundle $\varphi$ with a fiber $F$, then $\varphi_{|F|}$ coincides with $\varphi$ tautologically, so $[F]$ is sent back to $\varphi$.\fim
\newpage
\begin{coro}\label{properties_of_conic_bundle_fibers}
Every fiber of a conic bundle $\varphi:X\to\P^1$ is connected. Moreover, all fibers $\varphi$ are smooth, irreducible curves of genus $0$ except for finitely many which are reducible and supported on negative curves.
\end{coro}
\noindent\textit{Proof.} Let $F$ be any fiber of $\varphi$. By Theorem \ref{correspondence}, $[F]$ is a conic class and $\varphi_{|F|}:X\to\P^1$ is
precisely $\varphi$. Then $F$ is connected by Proposition \ref{connectedness} and the rest follows from Proposition \ref{finitely_many_sinuglar_fibers}.\fim
\\ \\
\indent Corollary \ref{properties_of_conic_bundle_fibers} is the starting point for the classification of conic bundle fibers in Section \ref{Classification}.

\section{Classification of conic bundle fibers}\label{Classification}\
\indent Let $\pi:X\to\P^1$ be a rational elliptic fibration. In this section we give a complete description of fibers of a conic bundle based on Corollary \ref{properties_of_conic_bundle_fibers} in Section \ref{NumChar}.

\begin{lemma}\label{lemma_structure}
Let $\varphi:X\to\P^1$ be a conic bundle and $D$ any fiber of $\varphi$. Then $D$ is connected and 
\begin{enumerate}[(i)]
\item $D$ is either a smooth, irreducible curve of genus $0$, or
\item $D=P_1+P_2+D'$, where $P_1,P_2$ are $(-1)$-curves (sections of $\pi$), not necessarily distinct, and $D'$ is either zero or supported on $(-2)$-curves (components of reducible fibers of $\pi$).
\end{enumerate}
\end{lemma}
\noindent\textit{Proof.} By Corollary \ref{properties_of_conic_bundle_fibers}, all fibers of $\varphi$ fall into category (i) except for finitely many that are reducible and supported on negative curves. Let $D$ be one of such finitely many.

From Lemma \ref{neg_curves}, $\text{Supp }D$ has only $(-1)$-curves (sections of $\pi$) or $(-2)$-curves (components of reducible fibers of $\pi$). By adjunction, if $C\in\text{Supp }D$ is a $(-2)$-curve, then $C\cdot(-K_X)=0$, and if $P\in\text{Supp }D$ is a $(-1)$-curve, then $P\cdot(-K_X)=1$. But $D\cdot(-K_X)=2$, so we must have $D=P_1+P_2+D'$ with the desired composition.\fim
\\ \\
\indent At this point we have enough information about the curves that support a conic bundle fiber. It remains to investigate their multiplicities and how they intersect one another.
\newpage

\begin{teor}\label{teor_classif}
Let $X$ be a rational elliptic surface with elliptic fibration $\pi:X\to\P^1$ and let $\varphi:X\to\P^1$ be a conic bundle. If $D$ is a fiber of $\varphi$, then the intersection graph of $D$ fits one of the types below. Conversely, if the intersection graph of a divisor $D$ fits any of these types, then $|D|$ induces a conic bundle $\varphi_{|D|}:X\to\P^1$.
\begin{table}[h]
\begin{center}
\centering
\begin{tabular}{|c|c|} 
\hline
Type & Intersection Graph\\ %& $D$\\
\hline
\multirow{3}{*}{\hfil $0$} & \multirow{3}{*}{\hfil \TypeE}\\ %& \multirow{3}{*}{\hfil $C$}\\
&\\
&\\
\hline
\multirow{3}{*}{\hfil $A_2$} & \multirow{3}{*}{\hfil \TypeD}\\%& \multirow{3}{*}{\hfil $P+Q$}\\
& \\
& \\
\hline
\multirow{3}{*}{\hfil $A_n$ ($n\geq 3$)} & \multirow{3}{*}{\hfil \TypeA}\\ %& \multirow{3}{*}{\hfil $P+\sum_i\Theta_i+Q$}\\
& \\ 
& \\
\hline
\multirow{3}{*}{\hfil $D_3$} & \multirow{3}{*}{\hfil \TypeC}\\ %& \multirow{3}{*}{\hfil $\Theta_0^1+2P+\Theta_0^2$}\\ 
& \\
& \\
\hline
\multirow{4}{*}{\hfil $D_m$ ($m\geq 4$)} & \multirow{4}{*}{\hfil \TypeB}\\ %& \multirow{4}{*}{\hfil $2P+\sum_{i=0}^{n-2}2\Theta_i+\Theta_{n-1}+\Theta_n$}\\
& \\ 
& \\
& \\
\hline
\end{tabular}
\end{center}
\begin{align*}
\star&\,\,\text{smooth, irreducible curve of genus }0\\
\circ&\,(-1)\text{-curve (section of }\pi)\\
\bullet&\,(-2)\text{-curve (component of a reducible fiber of }\pi)
\end{align*}
\end{table}
\end{teor}

\noindent\textbf{Terminology.}\label{terminology} Before we prove Theorem \ref{teor_classif}, we introduce a natural terminology for dealing with the intersection graph of $D$. When $C,C'\in\text{Supp }D$ are distinct, we say that $C'$ is a \textit{neighbour} of $C$ when $C\cdot C'>0$. If $C$ has exactly one neighbour, we call $C$ an \textit{extremity}. We denote the number of neighbours of $C$ by $n(C)$. A simple consequence of these definitions is the following lemma.

\begin{lemma}\label{neighbours}
If $D=\sum_in_iE_i$ is a fiber of a morphism $X\to\P^1$, then $n(E_i)\leq -n_i^2E_i^2$ for every $i$.
\end{lemma}
\noindent\textit{Proof.} By definition of $n(E_i)$, clearly $n(E_i)\leq \sum_{j\neq i}E_i\cdot E_j$. Since $D\cdot E_i=0$ by Lemma~\ref{lemma_fibers}, then
\begin{align*}
0=D\cdot E_i&=\sum_jn_jE_j\cdot E_i\\
&=n_iE_i^2+\sum_{j\neq i}n_jE_j\cdot E_i\\
&\geq n_iE_i^2+\sum_{j\neq i}E_j\cdot E_i\\
&\geq n_iE_i^2+n(E_i).\text{\fim}
\end{align*}
\newpage
\noindent\textit{Proof of Theorem \ref{teor_classif}.} We begin by the converse. If $D$ fits one of the types, we must prove $[D]\in\text{NS}(X)$ is a conic class, so that $\varphi_{|D|}:X\to\P^1$ is a conic bundle by Theorem \ref{correspondence}. A case-by-case verification gives $D\cdot C=0$ for all $C\in\text{Supp }D$. Since $D$ is effective, it is nef with $D^2=0$ by Lemma~\ref{critical_nef}. The condition $D\cdot(-K_X)=2$ is satified in type $0$ by adjunction. We have $D\cdot(-K_X)=2$ in types $A_2,A_n$ for they contain two distinct sections of $\pi$ and also in types $D_3,D_m$ for they contain a section with multiplicity 2. Hence $[D]\in\text{NS}(X)$ is a conic class, as desired.

Now let $D$ be a fiber of $\varphi$. By Lemma~\ref{lemma_structure}, $D$ is connected and has two possible forms. If $D$ is irreducible, we get type $0$. Otherwise $D=P_1+P_2+D'$, where $P_1,P_2$ are $(-1)$-curves and $D'$ is either zero or supported on $(-2)$-curves. If $D'=\sum_in_iC_i$ then $n(C_i)\leq 2n_i$ by Lemma~\ref{neighbours}. The bounds for $n(P_1), n(P_2)$ depend on whether i) $P_1\neq P_2$ or ii) $P_1=P_2$. In what follows Lemma~\ref{lemma_fibers}~a) is used implicitly several times.

\indent i) $P_1\neq P_2$. In this case $n(P_1)\leq 1$ and $n(P_2)\leq 1$. Since $D$ is connected, both $P_1,P_2$ must have some neighbour, so $n(P_1)=n(P_2)=1$, therefore $P_1,P_2$ are extremities. If the extremities $P_1,P_2$ meet, they form the whole graph, so $D=P_1+P_2$. This is type $A_2$.

If $P_1,P_2$ do not meet, by connectedness there must be a path on the intersection graph joining them, say $P_1,C_1,...,C_k,P_2$. Since $P_1$ is an extremity, it has only $C_1$ as a neighbour, so $0=D\cdot P_1=-1+n_1$ gives $n_1=1$. Moreover $n(C_1)\leq 2$ and by the position of $C_1$ in the path we have $n(C_1)=2$. We prove by induction that $n_i=1$ and $n(C_i)=2$ for all $i=1,...,k$. Assume this is true for $i=1,...,\ell<k$. Then $0=D\cdot C_\ell=1-2+n_{\ell+1}$, so $n_{\ell+1}=1$. Moreover, $n(C_{\ell+1})\leq 2$ and by the position of $C_{\ell+1}$ in the path we have $n(C_{\ell+1})=2$. So the graph is precisely the chain $P_1,C_1,...,C_k,P_2$. This is type $A_n$ ($n\geq 3$).
\\ \\
\indent ii) $P_1=P_2$. In this case $n(P_1)\leq 2$. We cannot have $n(P_1)=0$, otherwise $D^2=(2P_1)^2=-4$, so $n(P_1)=1$ or $2$. If $P_1$ has two neighbours, say $C_1,C_2$, then $0=D\cdot P_1=-2+n_1+n_2$, which only happens if $n_1=n_2=1$. Moreover, $n(C_1)\leq 2$, so $C_1$ can possibly have another neighbour $C_3$ in addition to $P_1$. But then $D\cdot P_1=0$ gives $n_3=0$, which is absurd, so $C_1$ has only $P_1$ as a neighbour. By symmetry $C_2$ also has only $P_1$ as a neighbour. This is type $D_3$.

Finally let $n(P_1)=1$ and $C_1$ be the only neighbour of $P_1$. Then $C_i\cdot P_1=0$ when $i>1$. Notice that $C_1,C_i$ come from the same fiber of $\pi$, say $F$, otherwise $C_i$ would be in a different connected component as $C_1,P_1$, which contradicts $D$ being connected. The possible Dynkin diagrams for $F$ are listed in Theorem \ref{Kodaira_fibers}. Since $P_1$ intersects $F$ in a simple component, the possibilities are

\begin{table}[h]
\begin{center}
\centering
\begin{tabular}{ccc} 
\multirow{4}{*}{\DiagramX} & \multirow{4}{*}{\hfil\DiagramY} & \multirow{4}{*}{\hfil\DiagramZ}\\
& &\\
& & 
\end{tabular}
\end{center}
\end{table}

In all diagrams above, $D\cdot P_1=0$ gives $n_1=2$. In the third diagram, $D\cdot C_1=0$ gives $n_2+n_m=2$, which only happens if $n_2=n_m=1$. But $D\cdot C_2=0$ gives $n_3=0$, which is absurd, so we rule out the third diagram. For the first two diagrams we proceed by induction: if $k<m$ and $n_1=...=n_k=2$, then $D\cdot C_k=0$ gives $n_{k+1}=2$, therefore $n_1=...=n_m=2$. But in the first diagram $D\cdot C_m=0$ gives $n_m=1$, which is absurd, so the first diagram is also ruled out.

Now let $C_{m+1},C_{m+2}$ be the first elements in branches $1$ and $2$ respectively. Then $D\cdot C_m=0$ gives $n_{m+1}=n_{m+2}=1$. Consequently $n(C_{m+1})\leq 2$ and $n(C_{m+2})\leq 2$. If $C_{m+1}$ has another neighbour $C_{m+3}$ in addition to $C_m$, then $D\cdot C_{m+1}=0$ gives $n_{m+3}=0$, which is absurd, so $C_{m+1}$ is an extremity. By symmetry, $C_{m+2}$ is also an extremity. This is type $D_m$ ($m\geq 4$).\fim
\newpage
\section{Fibers of conic bundles vs. fibers of the elliptic fibration}\label{conic_vs_elliptic}\
\indent Let $X$ be a rational elliptic surface with elliptic fibration $\pi:X\to\P^1$. The existence of a conic bundle $\varphi:X\to\P^1$ with a given fiber type is strongly dependent on the fiber configuration of $\pi$. This relationship is explored in Theorem~\ref{fiber_table}, which provides simple criteria to identify when a certain fiber type is possible. Before we prove it, we need the following result about the existence of disjoint sections.
\begin{lemma}\label{disjoint_sections}
If $X$ is a rational elliptic surface with at least two sections, then there exists a pair of disjoint sections.
\end{lemma}
\noindent\textit{Proof.} Let $\text{MW}(\pi)$ be the Mordell-Weil group of $\pi$, whose neutral section we denote by $O$. By \cite[Thm. 2.5]{OguisoShioda}, $\text{MW}(\pi)$ is generated by sections which are disjoint from $O$. Then there must be a generator $P\neq O$ disjoint from $O$, otherwise $\text{MW}(\pi)=\{O\}$, which contradicts the hypothesis.\fim
\\ \\
\indent We now state and prove the main result of this section.

\begin{teor}\label{fiber_table}
Let $X$ be a rational elliptic fibration with elliptic fibration $\pi:X\to\P^1$. Then the following statements hold:
\begin{enumerate}[a)]
\item $X$ admits a conic bundle with an $A_2$ fiber $\Rightarrow $ $\pi$ is not extremal (i.e. has positive generic rank).
\item $X$ admits a conic bundle with an $A_n$ ($n\geq 3$) fiber $\Leftrightarrow$ $\pi$ has a reducible fiber distinct from {\normalfont $\text{II}^*$}.
\item $X$ admits a conic bundle with a $D_3$ fiber $\Leftrightarrow$ $\pi$ has at least two reducible fibers.	
\item $X$ admits a conic bundle with a $D_m$ ($m\geq 4$) fiber $\Leftrightarrow$ $\pi$ has a nonreduced fiber or a fiber {\normalfont $\text{I}_n$} ($n\geq 4$).
\end{enumerate}
\end{teor}
\noindent\textit{Proof.} a) If $\pi$ is extremal, then all sections of $\pi$ are torsion sections, therefore no two of them ever meet \cite[Corollary 8.30]{MWL}, which makes type $A_2$ impossible.

b) Assume $X$ admits a conic bundle with an $A_n$ ($n\geq 3$) fiber. Since type $A_n$ contains a $(-2)$-curve, by Lemma \ref{neg_curves}, $\pi$ has a reducible fiber $F$. We claim that $F$ is not of type $\text{II}^*$. Indeed, if this were the case, the Mordell-Weil group $\text{MW}(\pi)$ would be trivial \cite{Persson}, which is impossible, since the $A_n$ fiber of the conic bundle contains two distinct sections. Conversely, assume $\pi$ has a reducible fiber $F$ other than $\text{II}^*$. Then $\text{MW}(\pi)$ is not trivial \cite{Persson} and by Lemma~\ref{disjoint_sections} we may find two disjoint sections $P,P'$. Let $C,C'\in\text{Supp }F$ be the components hit by $P,P'$. Since $F$ is connected, there is a path $C,C_1,...,C_\ell,C'$ in the intersection graph of $F$. Let $D:=P+C+C_1+...+C_\ell+C'+P'$. By Theorem \ref{teor_classif}, $\varphi_{|D|}:X\to\P^1$ is a conic bundle and $D$ is an $A_n$ fiber of it.
\begin{center}
\includegraphics[scale=0.65]{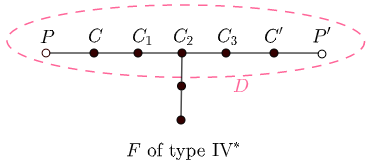}
\end{center}
\newpage
c) Assume $X$ admits a conic bundle with a $D_3$ fiber $D=C_1+2P+C_2$, where $C_1,C_2$ are $(-2)$-curves and $P$ is a section with $C_1\cdot P=C_2\cdot P=1$. Since $P$ hits each fiber of $\pi$ at exactly one point, then $C_1,C_2$ must come from two distinct reducible fibers of $\pi$. Conversely, let $F_1,F_2$ be two reducible fibers of $\pi$. If $P$ is a section, then $P$ hits $F_i$ at some $(-2)$-curve $C_i\in\text{Supp }F_i$. Let $D:=C_1+2P+C_2$. By Theorem \ref{teor_classif}, $\varphi_{|D|}:X\to\P^1$ is a conic bundle and $D$ is a $D_3$ fiber of it.
\begin{center}
\includegraphics[scale=0.65]{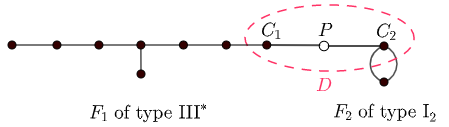}
\end{center}

d) Assume $X$ admits a conic bundle with a $D_m$ fiber $D=2P+2C_1+...+2C_\ell+(C_{\ell+1}+C_{\ell+2})$, where all $C_i$'s come from a reducible fiber $F$ of $\pi$. Notice that if $\ell> 1$ then $C_\ell$ meets three $(-2)$-curves, namely $C_{\ell-1},C_{\ell+1},C_{\ell+2}$ (see picture below). Going through the list in Theorem \ref{Kodaira_fibers}, we see that this intersection behavior only happens if $F$ is $\text{I}_n^*$, $\text{II}^*$, $\text{III}^*$ or $\text{IV}^*$, all of which are nonreduced. If $\ell=1$, then $C_1$ meets the section $P$ and two $(-2)$-curves which do not intersect, namely $C_2,C_3$. Again by examining the list in Theorem \ref{Kodaira_fibers}, $F$ must be $\text{I}_n$ with $n\geq 4$. Conversely, let $F$ be nonreduced or of type $\text{I}_n$ with $n\geq 4$. Take a section $P$ that hits $F$ at $C_1$. Now take a chain $C_2,...,C_\ell$ so that $C_\ell$ meets two other components of $F$. We name these two $C_{\ell+1},C_{\ell+2}$ and define $D:=2P+2C_1+...+2C_\ell+(C_{\ell+1}+C_{\ell+2})$. By Theorem \ref{teor_classif}, $\varphi_{|D|}:X\to\P^1$ is a conic bundle and $D$ is a type $D_m$ fiber of it.\fim
\begin{center}
\includegraphics[scale=0.65]{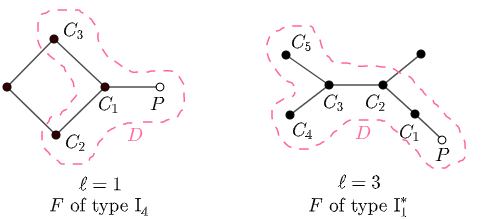}
\end{center}
\begin{remark}\label{converse_a)}
\normalfont The converse of Theorem \ref{fiber_table} a) is not true. Indeed, let $\pi$ with fiber configuration $(\text{III}^*,3\text{I}_1)$, in which case the generic rank is $1$. If $X$ admitted a conic bundle with an $A_2$ fiber, there would be sections $P_1,P_2$ with $P_1\cdot P_2=1$. However, by analyzing the Mordell-Weil lattice of $X$ and using the height formula \cite[Thm. 6.24]{MWL}, one finds that $P_1\cdot P_2=\frac{n^2-4}{4}$ or $\frac{n^2-1}{4}$ for some $n\in\Bbb{Z}$. In both cases $P_1\cdot P_2=1$ is impossible. A further study on the possible intersection numbers of sections is found in \cite{Costa}.
\end{remark}
\newpage
\section{Construction of conic bundles from a pencil of genus zero curves}\label{construction}\
\indent Let $X$ be a rational elliptic surface with elliptic fibration $\pi:X\to\P^1$. As described in Subsection~\ref{RES}, $\pi$ is induced by a pencil of cubics $\mathcal{P}$ from the blowup $p:X\to\P^2$ of the base locus of $\mathcal{P}$. We describe a method for constructing a conic bundle $\varphi:X\to\P^1$ from a pencil of curves with genus zero on $\P^2$. 

\noindent\textbf{Construction.} Let $\mathcal{Q}$ be a pencil of conics (or a pencil of lines) given by a dominant rational map $\psi:\P^2\dasharrow\P^1$ with the following properties:
\begin{enumerate}[(a)]
\item $\psi^{-1}(t)$ is smooth for all but finitely many $t\in\P^1$, i.e. the general member of $\mathcal{Q}$ is smooth. 
\item The indeterminacy locus of $\psi$ (equivalently, the base locus of $\mathcal{Q})$ is contained in the base locus of $\mathcal{P}$, which includes infinitely near points.
\end{enumerate}

Now define a surjective morphism $\varphi:X\to\P^1$ by the composition 
\begin{displaymath}
\begin{tikzcd}
X\arrow{r}{p}\arrow[bend right, swap]{rr}{\varphi} & \P^2\arrow[dashed]{r}{\psi} &\P^1
\end{tikzcd}
\end{displaymath}

Notice that $\varphi$ is a well defined conic bundle. Indeed, by property (b) the points of indeterminacy of $\psi$ are blown up under $p$, so $\varphi$ is a morphism. By property (a) the general fiber of $\varphi$ is a smooth, irreducible curve. Since $\mathcal{Q}$ is composed of conics (or lines), the general fiber of $\varphi$ has genus zero.
\\ \\
\noindent\textbf{Illustration.} Let $C$ be a smooth cubic and let $L_1,L_2,L_3$ be concurrent lines. Define $\mathcal{P}$ as the pencil generated by $C$ and $L_1+L_2+L_3$. Take points $P_1,P_2\in L_1\cap C$ and $P_3,P_4\in L_3\cap C$ and let $\mathcal{Q}$ be pencil of conics through $P_1,P_2,P_3,P_4$. In the following picture, $Q\in\mathcal{Q}$ is a general conic, so the strict transform of $Q$ under $p$ is a general fiber of the conic bundle $\varphi:X\to\P^1$.
\begin{center}
\includegraphics[scale=0.5]{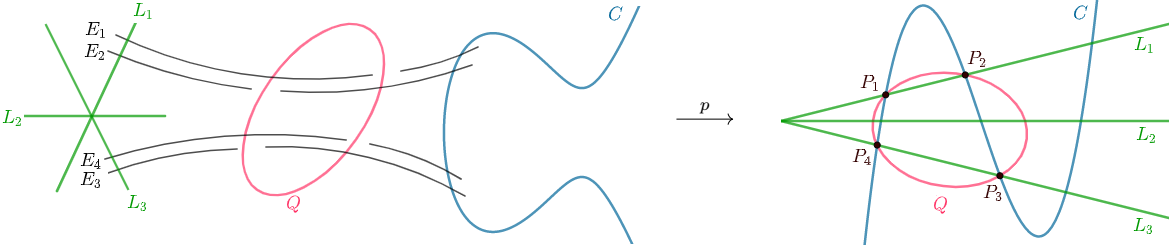}
\end{center}

\begin{remark}
\normalfont Since the base points $P_1,P_2,P_3,P_4$ of $\mathcal{Q}$ are blown up under $p$, then the pullback pencil $p^*\mathcal{Q}$ has four fixed components, namely the exceptional divisors $E_1,E_2,E_3,E_4$. By eliminating these we obtain a base point free pencil $p^*\mathcal{Q}-E_1-E_2-E_3-E_4$, which is precisely the one given by $\varphi:X\to\P^1$.
\end{remark}
\newpage
\section{Examples}\label{examples}\
\indent Let $\pi:X\to\P^1$ be a rational elliptic fibration. For simplicity we assume $\text{char}(k)=0$ throughout Section \ref{examples}, although similar constructions are possible in positive characteristic. We construct examples of conic bundles using the method described in Section \ref{construction}. 
\\ \\
\noindent\textbf{Notation.}
\begin{align*}
\mathcal{P}&\,\,\,\,\,\,\,\text{pencil of cubics on }\P^2\text{ inducing }\pi.\\
p&\,\,\,\,\,\,\,\text{blowup }p:X\to\P^2\text{ at the base locus of }\mathcal{P}.\\
\mathcal{Q}&\,\,\,\,\,\,\,\text{pencil of conic (or lines) on }\P^2.\\
Q&\,\,\,\,\,\,\,\text{conic in }\mathcal{Q},\text{ not necessarily smooth}.\\
L&\,\,\,\,\,\,\,\text{line in }\mathcal{Q}.\\
\varphi&\,\,\,\,\,\,\,\text{conic bundle }\varphi:X\to\P^1\text{ induced by }\mathcal{Q}\text{ (see Section \ref{construction})}.\\
D&\,\,\,\,\,\,\,\text{singular fiber of }\varphi\text{ such that }\varphi(D)=Q\text{ (or }\varphi(D)=L).\\
D'&\,\,\,\,\,\,\,\text{another singular fiber of }\varphi.
\end{align*}

\begin{remark}
\normalfont Since $D$ is a fiber of $\varphi$, by the correspondence in Proposition \ref{correspondence} the morphism $\varphi_{|D|}:X\to\P^1$ induced by $|D|$ coincides with $\varphi$ itself.
\end{remark}

\begin{remark}
\normalfont To simplify our notation, the strict transform of a curve $E\subset \P^2$ under $p$ is also denoted by $E$ instead of the usual $\widetilde{E}$.
\end{remark}

\subsection{Extreme cases}\
\indent By Theorem \ref{fiber_table}, there are two extreme cases in which $X$ can only admit conic bundles with exactly one type of singular fiber.
\begin{enumerate}[(1)]
\item When $\pi$ has a $\text{II}^*$ fiber: $X$ only admits conic bundles with singular fibers of type $D_m$ ($m\geq 4$).
\item When $\pi$ has no reducible fibers: $X$ only admits conic bundles with singular fibers of type $A_2$.
\end{enumerate}

In Persson's classification \cite{Persson}, case (1) corresponds to the first two entries in the list (the ones with trivial Mordell-Weil group) and case (2) corresponds to the last six (the ones with maximal Mordell-Weil rank). Examples \ref{only_type_5} and \ref{only_type_2} illustrate cases (1) and (2) respectively.
\newpage
\begin{example}\label{only_type_5}
We consider $\pi$ with configuration $(\text{II}^*,\text{II})$. By Theorem~\ref{fiber_table}, $X$ can only admit conic bundles with singular fibers of type $D_m$ ($m\geq 4$). We construct the elliptic surface $X$ by blowing up the base locus $\{9P_1\}$ of the pencil of cubics $\mathcal{P}$ induced by a smooth cubic {\color{customBlue}$C$} and a triple line {\color{customGreen}$3L$}. Let $\mathcal{Q}$ be the pencil of lines through $P_1$. Then ${\color{customPink}D}:=p^*{\color{customPink}L}-E_1$ is a curve in the base point free pencil $p^*\mathcal{Q}-E_1$, which induces the conic bundle $\varphi_{|{\color{customPink}D}|}:X\to\P^1$. The curve ${\color{customPink}D}$ is a $D_9$ fiber of $\varphi_{|{\color{customPink}D}|}$.
\begin{center}
\includegraphics[scale=0.65]{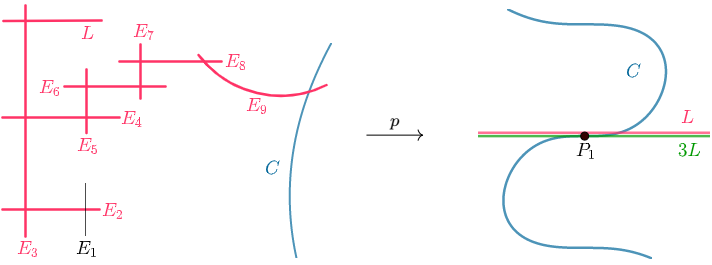}
\end{center}
\begin{align*}
&{\color{customBlue}C}:x^3+y^3-yz^2=0.\\
&{\color{customPink}L}:y=0.\\
&\mathcal{P}:\text{pencil of cubics induced by }{\color{customBlue}C}\text{ and }{\color{customGreen}3L}.\\
&\mathcal{Q}:\text{pencil of lines through }P_1.\\
&{\color{customPink}D}=2E_9+2E_8+2E_7+2E_6+2E_5+2E_4+2E_3+(E_2+{\color{customPink}L}), \text{ type }D_9.\\
&\text{Sequence of contractions:}\,E_9,E_8,E_7,E_6,E_5,E_4,E_3,E_2,E_1.
\end{align*}
\end{example}

\begin{remark}\label{single_singular_fiber}
\normalfont The conic bundle in Example~\ref{only_type_5} is in fact the only conic bundle on $X$. This follows from the fact that $\text{II}^*$ is the only reducible fiber of $\pi$ and that $E_9$ is the only section of $X$ since the Mordell-Weil group $\text{MW}(\pi)$ is trivial \cite{Persson}. By examining the intersection graph of $\text{II}^*$, ${\color{customPink}D}$ is the only divisor on $X$ forming a fiber of type $D_m$ ($m\geq 4$).
\end{remark}

\begin{remark}
\normalfont The case when $\pi$ has configuration $(\text{II}^*,2\text{I}_1)$ admits a similar construction and Remark \ref{single_singular_fiber} also applies.
\end{remark}
\newpage
\begin{example}\label{only_type_2}
We take $\pi$ with configuration $(\text{II},10\text{I}_1)$. By Theorem \ref{fiber_table}, $X$ can only admit conic bundles with singular fibers of type $A_2$. We construct the elliptic surface $X$ by blowing up the base locus $\{P_1,...,P_9\}$ of the pencil of cubics $\mathcal{P}$ induced by {\color{customBlue} $C$} and {\color{customGreen}$C'$}. Let $\mathcal{Q}$ be the pencil of lines through $P_1$. Then {\color{customPink}$D$}$:=p^*{\color{customPink}L}-E_1$ is a curve in the base point free pencil $p^*\mathcal{Q}-E_1$, which induces the conic bundle $\varphi_{|{\color{customPink}D}|}:X\to\P^1$. The curve ${\color{customPink}D}$ is an $A_2$ fiber of $\varphi_{|{\color{customPink}D}|}$.
\begin{center}
\includegraphics[scale=0.65]{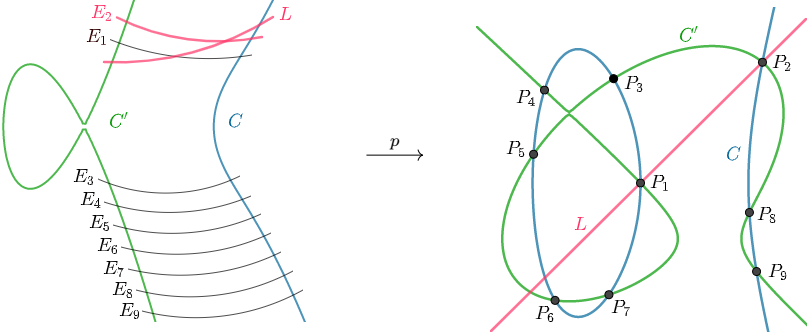}
\end{center}
\begin{align*}
&{\color{customBlue}C}:y^2z-4x^3+4xz^2=0.\\
&{\color{customGreen}C'}:y^2z-4x^3+4xz^2+(127/100)(xz^2-4y^3+4yz^2)=0.\\
&{\color{customPink}L}:\text{ line through }P_1,P_2.\\
&\mathcal{P}:\text{pencil of cubics induced by }{\color{customBlue}C}\text{ and }{\color{customGreen}C'}.\\
&\mathcal{Q}:\text{pencil of lines through }P_1.\\
&{\color{customPink}D}={\color{customPink}L}+E_2, \text{ type }A_2.\\
&\text{Sequence of contractions:}\,E_9,E_8,E_7,E_6,E_5,E_4,E_3,E_2,E_1.
\end{align*}
\end{example}

\begin{remark}
\normalfont ${\color{customPink}D}$ is not the only singular fiber of the conic bundle in Example~\ref{only_type_2}. In fact, each line $L_{1i}$ joining $P_1,P_i$ for any $i=2,...,9$ induces an $A_2$ fiber of $\varphi_{|{\color{customPink}D}|}$, namely $L_{1i}+E_i$.
\end{remark}

\begin{remark}
\normalfont Since the conic bundle in Example~\ref{only_type_2} only admits singular fibers of type $A_2$ and such fibers are isomorphic to a pair of lines meeting at a point, this is a \textit{standard conic bundle} in the sense of Manin and Tsfasman \cite[Subsection 2.2]{ManTsfa}.
\end{remark}
\newpage
\subsection{Mixed fiber types}
\begin{example}\label{types_2_and_3}
We take $\pi$ with configuration $(\text{IV}, \text{II}, 6\text{I}_1)$. By Theorem \ref{fiber_table}, $X$ admits only conic bundles with singular fibers of types $A_2$ or $A_n$ ($n\geq 3$). We construct a conic bundle with singular fibers of types $A_2$ and $A_3$. We construct the elliptic surface $X$ by blowing up the base locus $\{P_1,...,P_9\}$ of the pencil of cubics induced by {\color{customBlue}$C$} and {\color{customGreen}$L_1+L_2+L_3$}. Let $\mathcal{Q}$ be the pencil of lines through $P_1$. Then ${\color{customPink}D}:=p^*{\color{customPink}L}-E_1$ and ${\color{customPink}D'}:=p^*{\color{customPink}L'}-E_1$ are curves in the base point free pencil $p^*\mathcal{Q}-E_1$, which induces the conic bundle $\varphi_{|{\color{customPink}D}|}:X\to\P^1$. The curves $\color{customPink}D$, $\color{customPink}D'$  are fibers of $\varphi_{|\color{customPink}D|}$ of type $A_2$, $A_3$ respectively.
\begin{center}
\includegraphics[scale=0.65]{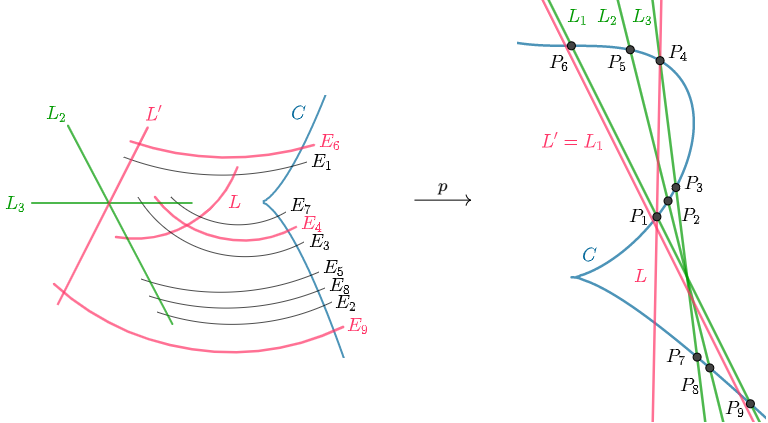}
\end{center}
\begin{align*}
&{\color{customBlue}C}:x^3+y^3-y^2z=0.\\
&{\color{customGreen}L_1}:y+2x-z=0.\\
&{\color{customGreen}L_2}:y+4x-2z=0.\\
&{\color{customGreen}L_3}:y+8x-4z=0.\\
&{\color{customPink}L}:\text{ line through }P_1,P_4.\\
&{\color{customPink}L'}={\color{customGreen}L_1}.\\
&\mathcal{P}:\text{ pencil of cubics induced by }{\color{customBlue}C}\text{ and }{\color{customGreen}L_1+L_2+L_3}.\\
&\mathcal{Q}:\text{ pencil of lines through }P_1.\\
&{\color{customPink}D}={\color{customPink}L}+E_4,\text{ type }A_2.\\
&{\color{customPink}D'}=E_6+{\color{customPink}L'}+E_9,\text{ type }A_3.\\
&\text{Sequence of contractions:}\,E_9,E_8,E_7,E_6,E_5,E_4,E_3,E_2,E_1.
\end{align*}
\end{example}
\begin{remark}
\normalfont In addition to ${\color{customPink}D}$ and ${\color{customPink}D'}$, the conic bundle in Example \ref{types_2_and_3} has five other singular fibers, each of type $A_2$. Namely, each line $L_{1i}$ joining $P_1,P_i$ with $i\in\{2,3,5,7,8\}$ induces the $A_2$ fiber $L_{1i}+E_i$.
\end{remark}

%\newpage
%\begin{example}
%Take $\pi$ with configuration $(\text{I}^*_2, 2\text{I}_2)$. By Theorem \ref{fiber_table}, $X$ only admits conic bundles with singular fibers of types 4 or 5. We construct a conic bundle with both fiber types.
%\begin{center}
%\includegraphics[scale=0.7]{example_8.2.6}
%\end{center}
%\begin{align*}
%&C:(x^2+y^2-1)\cdot y=0.\\
%&L_1:3x+4y+5=0.\\
%&L_2:3x-4y+5=0.\\
%&L:=L_1.\\
%&\mathcal{P}:\text{ induced by }C\text{ and }L_1+2L_2.\\
%&\mathcal{Q}:\text{ lines through }P_1.\\
%&D=E_2+2F_2+L_1,\text{ type 4}.\\
%&D'=2H_3+2G_3+2F_3+(E_3+L_2),\text{ type 5}.\\
%&\text{Contractions: }H_3,G_3,F_3,E_3,F_2,E_2,G_1,F_1,E_1.
%\end{align*}
%\end{example}
\newpage
%\subsection{Three types of singular fibers in the same conic bundle}
\begin{example}
We take $\pi$ with configuration $(\text{I}_7, \text{II}, 3\text{I}_1)$. By Theorem \ref{fiber_table}, $X$ admits conic bundles with singular fibers of types $A_n$ ($n\geq 3$), $D_m$ ($m\geq 4$) and possibly $A_2$. We construct a conic bundle with types $A_2$, $A_4$, $D_4$. The elliptic surface $X$ is constructed by blowing up the base locus $\{2P_1,3P_2, 2P_3,P_4,P_5\}$ of the pencil of cubics $\mathcal{P}$ induced by {\color{customBlue}$C$} and {\color{customGreen}$L_1+L_2+L_3$}. Let $\mathcal{Q}$ be the pencil of lines through $P_1$. Then ${\color{customPink}D}:=p^*{\color{customPink}L}-E_1$, ${\color{customPink}D'}:=p^*{\color{customPink}L'}-E_1$ and ${\color{customPink}D''}:=p^*{\color{customPink}L''}-E_1$ are curves in the base point free pencil $p^*\mathcal{Q}-E_1$, which induces the conic bundle $\varphi_{|{\color{customPink}D}|}:X\to\P^1$. The curves $\color{customPink}D$, $\color{customPink}D'$, $\color{customPink}D''$ are fibers of $\varphi_{|\color{customPink}D|}$ of types $D_4$, $A_4$, $A_2$ respectively.
\begin{center}
\includegraphics[scale=0.65]{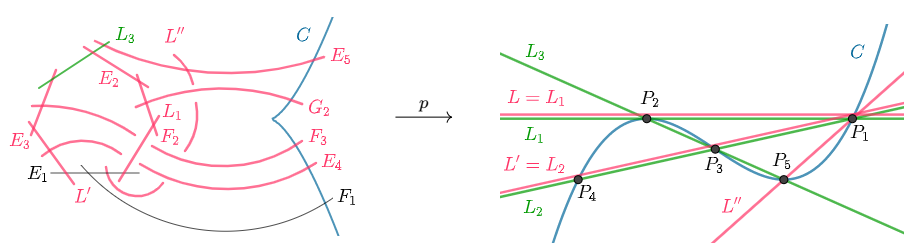}
\end{center}
\begin{align*}
&{\color{customBlue}C}:yz^2-2x^2(x-z)=0.\\
&{\color{customGreen}L_1}:y=0.\\
&{\color{customGreen}L_2}:2x-9y-2z=0.\\
&{\color{customGreen}L_3}:4x+9y=0.\\
&{\color{customPink}L}={\color{customGreen}L_1}.\\
&{\color{customPink}L'}={\color{customGreen}L_2}.\\
&{\color{customPink}L''}:\text{ line joining }P_1,P_5.\\
&\mathcal{P}:\text{ pencil of cubics induced by }{\color{customBlue}C}\text{ and }{\color{customGreen}L_1+L_2+L_3}.\\
&\mathcal{Q}:\text{ pencil of lines through }P_1.\\
&{\color{customPink}D}=2G_2+2F_2+(E_2+{\color{customPink}L}),\text{ type }D_4.\\
&{\color{customPink}D'}=E_4+{\color{customPink}L'}+E_3+F_3,\text{ type }A_4.\\
&{\color{customPink}D''}={\color{customPink}L''}+E_5,\text{ type }A_2.\\
&\text{Sequence of contractions:}\,E_5,E_4,F_3,E_3,G_2,F_2,E_2,F_1,E_1.
\end{align*}
\end{example}
\newpage
\begin{example}
We take $\pi$ with configuration $(\text{I}^*_2,\text{III},\text{I}_1)$. By Theorem \ref{fiber_table}, $X$ admits conic bundles
with singular fibers of types $A_n$ ($n\geq 3$), $D_3$, $D_m$ ($m\geq 4$) and possibly $A_2$. We construct a conic bundle with singular fiber of types $A_3$, $D_3$, $D_5$. The elliptic surface $X$ is constructed by blowing up the base locus $\{P_1,...,P_9\}$ of the pencil of cubics $\mathcal{P}$ induced by {\color{customBlue}$C$} and {\color{customGreen}$L_1+L_2+L_3$}. Let $\mathcal{Q}$ be the pencil of lines through $P_1$. Then ${\color{customPink}D}:=p^*{\color{customPink}L}-E_1$, ${\color{customPink}D'}:=p^*{\color{customPink}L'}-E_1$ and ${\color{customPink}D''}:=p^*{\color{customPink}L''}-E_1$ are curves in the base point free pencil $p^*\mathcal{Q}-E_1$, which induces the conic bundle $\varphi_{|{\color{customPink}D}|}:X\to\P^1$. The curves $\color{customPink}D$, $\color{customPink}D'$, $\color{customPink}D''$ are fibers of $\varphi_{|{\color{customPink}D}|}$ of types $A_3$, $D_3$, $D_5$ respectively.
\begin{center}
\includegraphics[scale=0.65]{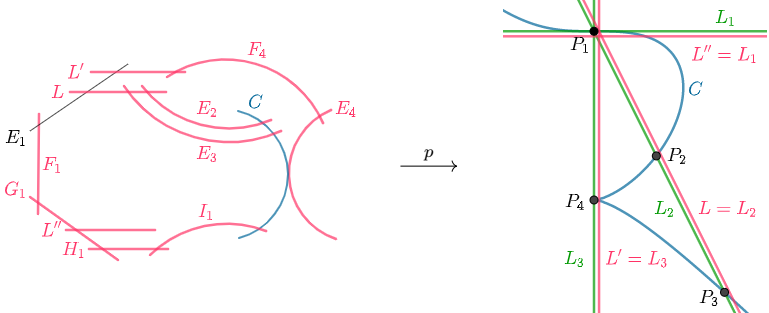}
\end{center}
\begin{align*}
&{\color{customBlue}C}:x^3+y^3-y^2z=0.\\
&{\color{customGreen}L_1}:y-z=0.\\
&{\color{customGreen}L_2}:x=0.\\
&{\color{customGreen}L_3}:y+2x-z	=0.\\
&{\color{customPink}L}:={\color{customGreen}L_2}.\\
&{\color{customPink}L'}:={\color{customGreen}L_3}.\\
&{\color{customPink}L''}:={\color{customGreen}L_1}.\\
&\mathcal{P}:\text{ pencil of cubics induced by }{\color{customBlue}C}\text{ and }{\color{customGreen}L_1+L_2+L_3}.\\
&\mathcal{Q}:\text{ pencil of lines through }P_1.\\
&{\color{customPink}D}=E_2+{\color{customPink}L}+E_3,\text{ type }A_3.\\
&{\color{customPink}D'}=E_4+2F_4+{\color{customPink}L'},\text{ type }D_3.\\
&{\color{customPink}D''}=2I_1+2H_1+2G_1+({\color{customPink}L''}+F_1),\text{ type }D_5.\\
&\text{Sequence of contractions:}\,F_4,E_4,E_3,E_2,I_1,H_1,G_1,F_1,E_1.
\end{align*}
\end{example}
%\newpage
%\subsection{Conic bundles with repeated types of singular fibers}
%\begin{example}
%Take $\pi$ with configuration $(\text{I}_8,4\text{I}_1)$. We construct a conic bundle with two singular fibers of type 3.
%\begin{center}
%\includegraphics[scale=0.68]{example_8.4.2}
%\end{center}
%\begin{align*}
%&C:y^2z-x^3+(1/36)x^2z+(5/36)xz^2-(25/1296)z^3=0.\\
%&L_1:126x-108y-35z=0.\\
%&L_2:18x-36y-5=0.\\
%&L_3:30x-36y-5=0.\\
%&\mathcal{P}:\text{ induced by }C\text{ and }L_1+L_2+L_3.\\
%&\mathcal{Q}:\text{ lines through }P_2.\\
%&D=G_2+F_2+L_1+E_1+F_1,\text{ type 3}.\\
%&D'=E_4+L_3+E_3+F_3+G_3,\text{ type 3}.\\
%&\text{Contractions: }E_4,G_3,F_3,E_3,G_2,F_2,E_2,F_1,E_1.
%\end{align*}
%\end{example}
%\newpage
%\begin{example}
%Take $\pi$ with configuration $(4\text{I}^*_0)$. We construct a conic bundle with four singular fibers of type 4.
%\begin{center}
%\includegraphics[scale=0.7]{example_8.4.3}
%\end{center}
%\begin{align*}
%&M:x=0.\\
%&N:y=0.\\
%&L_1:y+3x-z=0.\\
%&L_2:2y+3x-2z=0.\\
%&L_3:y+x-z=0.\\
%&\mathcal{P}:\text{ induced by }M+2N\text{ and }L_1+L_2+L_3.\\
%&\mathcal{Q}:\text{ lines through }P_1.\\
%&D=L_1+2F_2+E_2,\text{ type 4}.\\
%&D'=L_2+2F_3+E_3,\text{ type 4}.\\
%&D''=L_3+2F_4+E_4,\text{ type 4}.\\
%&D'''=M+2G_1+F_1,\text{ type 4}.\\
%&\text{Contractions: }F_4,E_4,F_3,E_3,F_2,E_2,G_1,F_1,E_1.
%\end{align*}
%\end{example}
\newpage
\bibliographystyle{alpha}
\bibliography{my}
\end{document}